\documentclass{amsart}
\usepackage{epsf, amssymb}

\newtheorem{prop}{Proposition}[section]
\newtheorem{theorem}{Theorem}
\newtheorem{lemma}[prop]{Lemma}

\newcommand{\R}{{\bf R}} 
 
\newcommand{\K}{{\mathcal K}} 
\newcommand{\s}{{\mathcal S}} 
\newcommand{\p}[1]{{\mathcal P}_{#1}} 
\newcommand{\hra}{\hookrightarrow} 

\title{On a Map From Pure Braids To Knots}

\author{Jacob Mostovoy}
\address{Instituto de Matem\'{a}ticas (Unidad Cuernavaca),
Universidad Nacional Aut\'{o}noma de M\'{e}xico,
A.P. 273-3 Admon. de correos no. 3, Cuernavaca, Morelos, MEXICO}
\email{jacob@matcuer.unam.mx}

\author{Theodore Stanford}
\address{Mathematics Department,
United States Naval Academy,
572C Holloway Road,
Annapolis MD 21402, USA}
\email{stanford@nadn.navy.mil}

\begin{document}

\vspace{2cm}
\maketitle

We study a certain type of braid closure which resembles the plat closure
but has certain advantages; for example, it maps pure braids to knots.
The main results of this note are a Markov-type theorem and a description of
how Vassiliev invariants behave under this braid closure.  


\section{Definition and properties of the short-circuit map.}

\begin{figure}[ht]
\[\epsffile{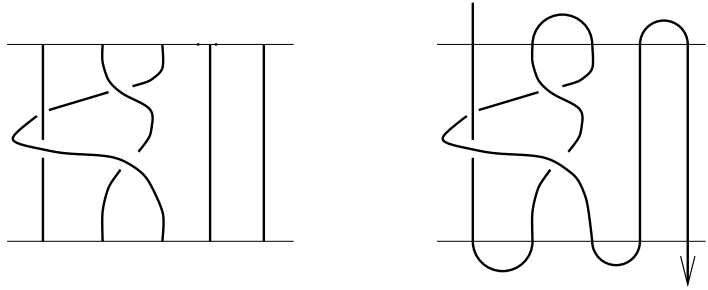}\]
\caption{A braid $b\in\p{5}$ and the knot $\s_2(b)$.}\label{f:def}
\end{figure}

We define the ``short-circuit'' map $\s_n$ from the pure
braid group on $2n+1$ strands $\p{2n+1}$ to the monoid of
the isotopy classes of oriented knots $\K$ as pictured on
Figure~\ref{f:def}. The strands of the braid are joined
together in turn at the bottom and at the top.
We think of knots as of non-compact, or ``long'' knots
here. These maps are compatible with the inclusions
$\p{2n+1}\hra\p{2n+3}$ so they extend to a map $\s
:\p{\infty}\to \K.$ Here by $\p{\infty}$ we understand the inductive limit 
of the sequence of inclusions $\p{i}\hra\p{i+1}$.
 
\begin{figure}[ht]
\[\epsffile{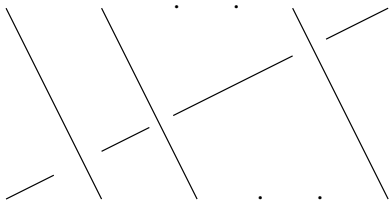}\]
\caption{The braid $t_n$.}\label{f:plat}
\end{figure}

The construction and, as we will see later, some properties
of the map $\s$ resemble those of the plat closure which
sends braids with even number of strands to links. (For the
definition and properties of the plat closure see
\cite{Bi1,Bi2}.) Indeed, if $t_n$ denotes the $2n$-strand
braid pictured on Figure~\ref{f:plat}, then for any
$x\in\p{2n+1}$ the (unoriented) knot $\s(x)$ is equivalent to the knot,
obtained by taking the image of $x$ in $\p{2n+2}$ under the
standard inclusion, multiplying by $t_{n+1}$ on the left
(i.e.\ on the top) and taking the plat closure.

However, if we are interested in knots rather than links the
map $\s$ is more convenient than the plat closure.
The most obvious difference is the behaviour under
stabilization maps and tensor products. Adding two
unbraided strands to a braid changes its image under the
plat closure by adding an unknotted and unlinked component,
while the image of the short-circuit map does not change. As
for tensor (external) products, the plat closure sends a
product of braids to the distant union of their plat
closures, while under short-circuiting the tensor product of
braids is sent to the connected sum of the corresponding
knots.

To make the last statement more precise, we may define the tensor
product of two pure braids with odd numbers of strands as follows.

Let $i:\p{2n+1}\hra\p{2(n+m)+1}$ be the standard inclusion onto the first 
$2n+1$ strands and $i':\p{2m+1}\hra\p{2(n+m)+1}$ be the inclusion onto
the last $2m+1$ strands.  Then we can define a product
\[ \p{2n+1}\otimes\p{2m+1}\to\p{2(n+m)+1} \]
by sending a pair $(b_1,b_2)$, where $b_1\in\p{2n+1}$ and
$b_2\in\p{2m+1}$ to
\[i(b_1)i'(b_2)\in\p{2(n+m)+1}.\] 
With this definition it is clear that 
\[ \s_n(b_1)\#\s_m(b_2)=\s_{n+m}(b_1\otimes b_2).\]

The restriction to an odd number of strands is 
by no means crucial.  If $b\in\p{2n}$ we can define an analogue of the 
short-circuit map as a suitably oriented plat closure of the braid $t_{n}b$. 
This definition is equally good for the purposes of our paper and has certain
advantages.  Namely, this version of the short-circuit
closure respects the usual tensor product of braids;
also, in this set-up Theorem~\ref{thm:bridge} below becomes
tautological.

Nevertheless, we prefer to work with braids on odd number of
strands.  It follows from Theorem~\ref{thm:bridge} that any
knot which can be realized as a plat closure of a
$2n$-stranded braid can be obtained by short-circuiting some
pure braid on $2n-1$ strands.  This generalizes the
well-known fact that a 2-bridge knot can be represented by a
braid in $\p{3}$. In this sense, the short-circuit map for
$\p{\rm odd}$ is more ``economic''.  We repeat, however,
that in our context this is a matter of taste.

\subsection{Filtration by the number of strands and the bridge number.}

Any filtration on the infinite pure braid group $\p{\infty}$
is sent by $\s$ to a filtration on knots. The most obvious
filtration on $\p{\infty}$ to consider is the filtration
``by the number of strands''
\[\p{1}\subset\p{3}\subset\p{5}\subset\ldots\subset\p{\infty}.\]
\begin{theorem}\label{thm:bridge}
The filtration on knots by $\s(\p{2n+1})$ is the filtration by knots
with bridge number less than or equal to $n+1$. 
\end{theorem}

To prove Theorem~\ref{thm:bridge} it is enough to show that the minimal
number of maxima of the height function in a realization of a knot in 
$\R^3$ as a long knot is the bridge number minus 1; this will be done in 
Section~\ref{s:bridge}.

The bridge number minus 1 is
an additive knot invariant (see \cite{Sch}) 
so, the filtration by
$\s(\p{2n+1})$ gives rise to an additive grading on $\K$.

\subsection{Structure of the short-circuit map.}

First we introduce some notation. By $A_{i,j}$ where $i\neq
j$ are positive integers we denote the standard generators
of $\p{\infty}$. By $\phi_i^{n}$ we mean the homomorphism
$\p{2n}\to\p{2n+1}$ which doubles the $i$th strand.
Homomorphisms $\phi_i^{n}$ respect the standard inclusions
of the pure braid groups so as $n$ tends to infinity the limit 
$\phi_i:\p{\infty}\to\p{\infty}$ 
is well-defined. 

\begin{figure}[hbt]
\[\epsffile{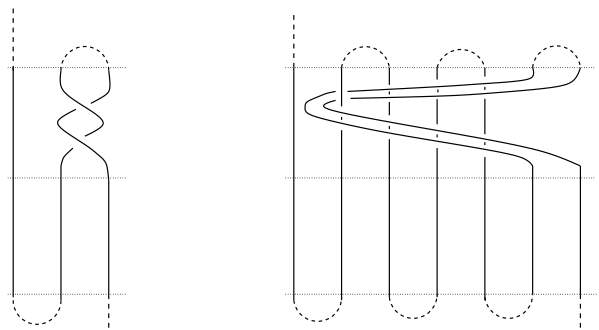}\]
\caption{
Examples of the action of $H^T$ on the trivial braid.}\label{f:act}
\end{figure}

Let $H^T\in\p{\infty}$ be the subgroup generated by
$A_{i,i+1}$ and $\phi_i(A_{i,j})$
for all even $i$ and all $j\neq i$.  Similarly we define the
subgroup $H^B$ with the only difference that $i$ is required
to be odd.
The subgroup $H^T$ acts on $\p{\infty}$ on the left and this action
preserves the fibres of $\s$, see Figure~\ref{f:act}. 
Similarly, $H^B$ act on $\p{\infty}$ on the right, also
preserving the fibres.

\begin{theorem}\label{thm:structure}
The short-circuit map identifies the monoid of knots $\K$ with the 
quotient set $H^T \backslash \p{\infty} /H^B$.
\end{theorem} 

This theorem is a version of the main theorem of \cite{Bi2}
which describes the equivalence classes of plat
closures. The proof we sketch in Section~\ref{s:structure}
is simplified by the fact the we are only interested in
knots.  Note also that Birman's theorem as stated in \cite{Bi2} concerns 
unoriented knot and link types, whereas our theorem concerns oriented
knot types.

\subsection{Lower central series and Vassiliev invariants.}

One can easily check that Vassiliev knot invariants pull
back under the short-circuit map to Vassiliev invariants of
braids. The action of $H^T$ and $H^B$ on $\p{\infty}$
induces an action on Vassiliev braid invariants which,
clearly, preserves the type. (Here we do not assume the
invariants to be normalized, i.e.\ do not require them to
take a prescribed value on the trivial braid.) Thus the
finite type knot invariants can be identified with those
finite type pure braid invariants which are fixed by the
two-sided action of $H^T$ and $H^B$.

Sometimes it is more convenient, however, to think of
Vassiliev invariants in the dual setting. Recall that a knot
(pure braid) is called $n$-trivial if it cannot be
distinguished from the the trivial knot (braid) by
invariants of order less than $n$. For pure braids
$n$-triviality is well-understood: $b\in\p{k}$ is
$n$-trivial if and only if $b\in\gamma_{n}\p{k}$ - the $n$-th
term of the lower central series of $\p{k}$.

Let $\K_n\subset\K$ be the set of $n$-trivial knots.
\begin{theorem}\label{thm:lcs}
Short-circuiting sends the filtration of $\p{\infty}$ by the
lower central series to the filtration by $n$-trivial knots:
\[\s(\gamma_{n}\p{\infty})=\K_n.\]
\end{theorem}

This allows to formulate problems from the theory of
Vassiliev knot invariants in purely group-theoretic
terms. For example, finite type knot invariants separate the
unknot if and only if any orbit of the two-sided action of
$H^T$ and $H^B$, apart from the orbit of the trivial braid,
intersects only a finite number of terms of the lower
central series.  Another way to state this is to consider
the nilpotent topology on $\p{\infty}$ (with basis 
the cosets of $\gamma_{n}{\p{\infty}}$ for all $n$).
Then finite type invariants separate the unknot if and
only if the set $H^T H^B = \{tb\ |\ t \in H^T, b \in H^B\}$
is closed in the nilpotent topology.

The proof of Theorem~\ref{thm:lcs} follows closely
the same arguments as in \cite{St}.
It is even simplified in some ways in our setting.
For example, if $x$ and $y$ are two braids, then
$\s{(x)} \# \s{(y)} = \s{(xtyb)} = \s{(xty)} = \s ((t^{-1}xtx^{-1})xy)$,
which is equivalent to $\s (xy)$ modulo a commutator.  Inductively,
braid product and connected sum are equivalent, modulo commutators
of higher order, which is the main idea behind the results in \cite{St}.


\section{Bridge number for long knots.}\label{s:bridge}

Here we will see that  the minimal
number of maxima $b_L$ 
of the ``height function'' in the realization of a knot in 
$\R^3$ as a long knot is less by 1 than  the minimal
number of maxima of the height function in the compact realization 
$S^1\hra \R^3$ of the same knot, i.e.\ than the bridge number $b$.

\begin{figure}[ht]
\[\epsffile{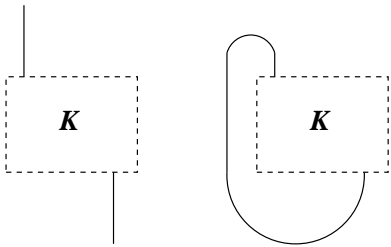}\]
\caption{A long knot $K$ and the corresponding compact knot.}\label{f:close}
\end{figure}

For a long knot with $b_L$ maxima of the height function it is obvious
that there exist a compact embedding of the same knot with  $b_L+1$
maxima, see Figure~\ref{f:close}.

Conversely, let $k$ be a compact knot $S^1\hra \R^3$ with $b$ maxima
and $b$ minima which can be taken to be non-degenerate. We construct a long 
knot $k'$ with  $b-1$ maxima which
is equivalent to $k$ as follows.

Choose a point on $k$ which is not critical for the height function to be the
origin in $\R^3$. Let $A$ be the maximum and $B$ the minimum between which the
chosen point  lies; by $AB$ we denote the closed segment of $k$ 
which lies between $A$ and $B$ and passes through the origin. 

Let $F(t):\R\to\R^3$ be a curve which intersects each horizontal plane once 
and such that its intersection with the knot $k$ is exactly the segment $AB$. 
We can assume that the curve $F$ is parametrized by the $z$-coordinate in 
$\R^3$, i.e. $F(t)=(F_x(t),F_y(t),t)$, and that $F$ is a smooth function of
$t$ everywhere apart from the points where $F(t)=A$ or $F(t)=B$.

\begin{figure}[ht]
\[\epsffile{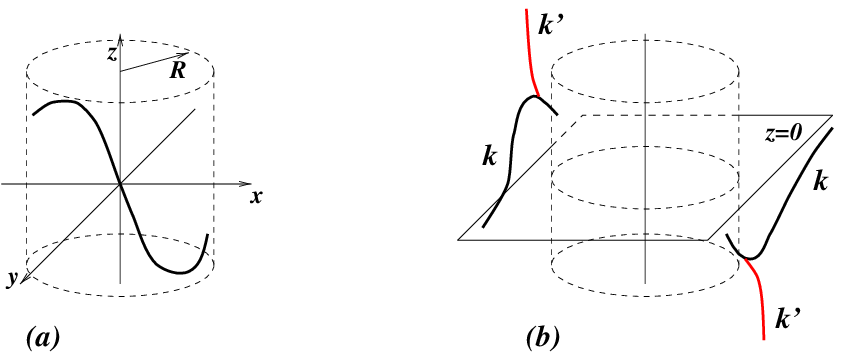}\]
\caption{}\label{f:cylinder}
\end{figure}

Consider a map $\Phi: \R^3\to\R^3$ given by
\[\Phi(x,y,z) = (x-F_x(z),y-F_y(z),z).\]
The transformation $\Phi$ preserves the horizontal planes, so it
does not change the number of maxima and minima of the height function on the
knot $k$. It is clear that there exist such $R>0$ that 
the intersection of the image of the embedding $\Phi(k)$ with
the cylinder $x^2+y^2<R^2$ is an interval, which is embedded 
with exactly one minimum and one maximum of the height function. Strictly 
speaking, the embedding $\Phi(k)$ is only piecewise-smooth, however, we
can smooth it out in such a way that its intersection with 
the cylindrical neighbourhood of the $z$-axis of radius $R$ is an interval,
which intersects the $z$-axis in the origin only and which is embedded 
with exactly one minimum and one maximum of the height function,
see Figure\ref{f:cylinder}(a). 

Thus in what follows we can assume that $k$
has the above form.
Now we compactify $\R^3$ to $S^3$ by an interval adding a point 
at infinity to each horizontal plane and  two points $z=\pm\infty$.
Denote by $V\subset S^3$ a copy of $\R^3$ obtained by throwing out the closure
of the $z$-axis. The intersection of the knot $k$ with $V$
is a long knot, which is equivalent to $k$ if we choose the 
orientation of $V$ to be compatible with that of $\R^3$. In the coordinates 
centred at the point at infinity whose $z$-coordinate is zero, this long
knot looks as on Figure~\ref{f:cylinder}(b). Obviously, it is 
equivalent to the knot $k'$ that differs from $k$ only inside the 
cylindrical neighbourhood of the $z$-axis (which 
is pictured as the outside part of the cylinder on Figure~\ref{f:cylinder}(b))
and has exactly $b-1$ maxima and $b-1$ minima.
  

\section{Short-circuit map as a two-sided quotient map.}\label{s:structure}

We say that a smooth long knot $k(t):\R\to\R^3$ is a Morse knot if the height 
function on it: (a) has only a finite number of critical points, all of which 
are non-degenerate; (b) tends to $\pm\infty$ as $t\to\mp\infty$; in other 
words, we assume that all knots ``point downwards''. Two Morse knots are Morse
equivalent if one can be deformed into the other through Morse knots.

\begin{figure}[ht]
\[\epsffile{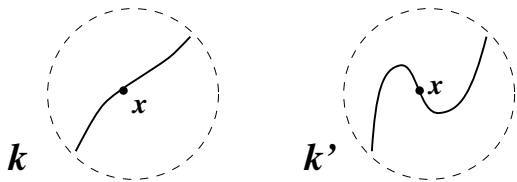}\]
\caption{Inserting a hump.}\label{f:insert}
\end{figure}

Let $k$ be a Morse knot and $x$ be a point on $k$ which is non-critical for
the height function. We will say that a knot $k'$ is obtained from $k$ by 
insertion of a hump at $x$ if $k$ and $k'$ coincide outside some small 
neighbourhood of $x$ and inside this neighbourhood they differ as on 
Figure~\ref{f:insert}.   

\begin{lemma}\label{lemma:insert}
Any two knots obtained from the same Morse knot by insertion of a hump 
are Morse equivalent. 
\end{lemma}
\begin{proof}
The lemma is clearly true if there are no critical points of the height 
function between the points $x_1$ and $x_2$ where we insert humps. 
In case there is one critical point between $x_1$ and $x_2$ the lemma follows
\begin{figure}[ht]
\[\epsffile{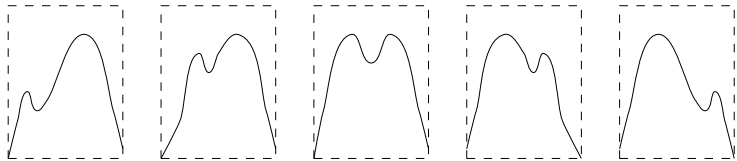}\]
\caption{Passing a hump through a critical point.}\label{f:movehump}
\end{figure}
from the argument on Figure~\ref{f:movehump}. This also proves the lemma
in the general case.
\end{proof}

Let $b_1\in\p{2n+1}$ and $b_2\in\p{2m+1}$ and, as before, denote by $i(b_k)$ 
the image
of the standard inclusion of $b_k$ into $\p{2N+1}$, $N\geqslant n,m$.

\begin{lemma}\label{lemma:Me}
If $\s_{n}(b_1)$ and $\s_{m}(b_2)$ are in the same isotopy class in $\K$ there
exists $N\geqslant n,m$ such that $\s_{N}(i(b_1))$ and $\s_{N}(i(b_2))$ are 
Morse equivalent.
\end{lemma}
\begin{proof}
Let 
\[ f^T(t)=(f_x^T(t),f_y^T(t),f_z^T(t))\] 
where $T\in [0,1]$ and $t\in\R$ be a homotopy between
$\s_{n}(b_1)$ and $\s_{m}(b_2)$, that is, for each $T$ the map 
$f^T(t):\R\to\R^3$
defines a long knot and $f^0(t)=\s_{n}(b_1)$ and $f^1(t)=\s_{m}(b_2)$.

In $[0,1]\times\R$ consider the subset $W$ of pairs $(T,t)$ such that 
$\frac{\partial}{\partial{t}} f_z^T(t)=0.$
Without loss of generality we can assume that $W$ is a union of smooth 
compact non-singular curves whose boundary is either empty or belongs to 
$(\{0\}\cup\{1\})\times\R$ and that there are only a finite number of 
tangencies of $W$ with horizontal lines of the form $\{T\}\times\R$. In
addition we require these tangencies to take place at different values
of the parameter $T$; see Figure~\ref{f:homotopy}. 
These assumptions imply, in particular, that for
all but a finite number of values of $T$ the knot $f^T(t)$ is Morse and
that the perestroikas at the bifurcation values of $T$ are generic, 
i.e.\ are insertions (or removals) of humps.

If there are no points of tangency of $W$ with  horizontal lines
the knots $\s_{n}(b_1)$ and $\s_{m}(b_2)$ are Morse equivalent and $n=m=N$.

Otherwise,
choose the point of tangency of $W$ with a horizontal line which corresponds 
to the insertion of a hump with the smallest value of $T$.
It is clear that we can connect it with the lower 
boundary line $\{0\}\times\R$ by a segment $s$ of a curve which is  
disjoint from $W$ and whose tangent is nowhere
horizontal, see Figure~\ref{f:homotopy}(a).
\begin{figure}[ht]
\[\epsffile{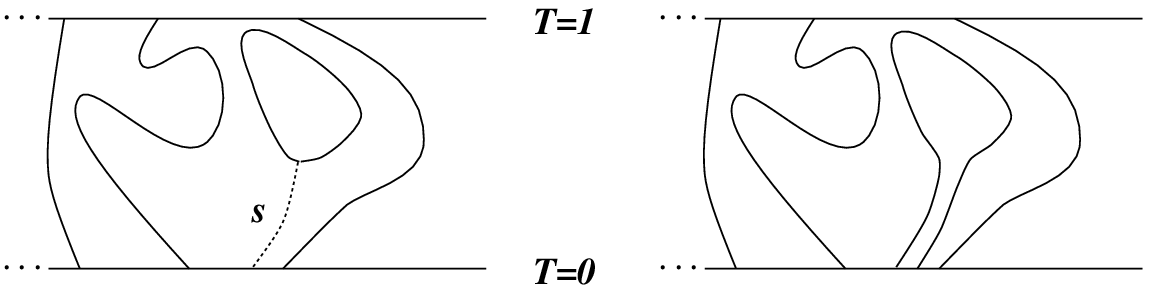}\]
\caption{}\label{f:homotopy}
\end{figure}

In the  neighbourhood of each point of $s$ we can modify the knots
$f^T(t)$ by inserting humps, this changes $W$ as shown on 
Figure~\ref{f:homotopy}(b). Notice that the number of points where $W$ has a
horizontal tangent has decreased by one and the knot $f^0(t)=\s_{n}(b_1)$
was changed by an insertion of a hump.

Thus, proceeding inductively, we eliminate all insertions of humps.
In the same way we eliminate the removals of humps with the only 
difference that we connect them to the upper boundary line and proceed
from the bifurcation with the largest value of $T$ downwards.

The result is that we construct a Morse equivalence between  $\s_{n}(b_1)$,
possibly with several humps inserted, and $\s_{m}(b_2)$, also with some 
extra humps. However, from Lemma~\ref{lemma:insert} we know that $\s_{n}(b_1)$
and $\s_{m}(b_2)$ with humps inserted are  Morse equivalent to 
$\s_{N}(i(b_1))$ and $\s_{N}(i(b_2))$ respectively (here $N$ is the number
of maxima of the modified knots) and this proves the lemma.

\end{proof}


Let $b_1\in\p{2N+1}$ and $b_2\in\p{2N+1}$ represent the same knot. 
Lemma~\ref{lemma:Me} allows us to assume that the knots $\s_N(b_1)$ and 
$\s_N(b_1)$ are Morse equivalent.

Given a deformation of $\s_N(b_1)$ to $\s_N(b_1)$ through Morse knots 
we are going to construct a one-dimensional family of 
braids $f^T:[0,1]\to\p{2N+1}$ such that $f^0=b_1$, $f^1=b_2$ and which is not
continuous only at a finite number of values of the parameter, where the 
``jump'' can be expressed as the multiplication by some element of $H^T$ or 
$H^B$.

\begin{figure}[ht]
\[\epsffile{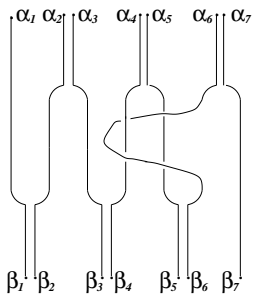}\]
\caption{Getting a braid from a knot.}\label{f:susp}
\end{figure}

The braid $f^0$ is obtained by ``suspending'' the knot $\s_N(b_1)$ by maxima 
and minima, see Figure~\ref{f:susp}. Here we choose the points $\alpha_i$
and $\beta_i$ in such a way that the deformation of $\s_n(b_1)$ into 
$\s_n(b_2)$ takes place entirely between the horizontal planes in which 
$\alpha_i$ and $\beta_i$ are situated. Of course, $f^0$ is the same braid as 
$b_1$. Think of the double lines which connect maxima and minima with the 
points $\alpha_i$ and $\beta_i$ respectively as of very narrow rubber strips. 
Then, if we deform the knot keeping  the points $\alpha_i$ and $\beta_i$
fixed,  the suspended knot also deforms and gives the 
braid $f^T$. 

It may happen in the process of deformation that some rubber strips 
intersect the knot or intersect each other. Without loss of generality we can 
assume that these events take place near a finite number of distinct 
values of $T$.

Suppose that the rubber strip which 
connects a maximum with points $\alpha_i$ and $\alpha_{i+1}$ intersects the 
knot between $T=T_{0}$ and $T=T_{0}+\epsilon$. Then one can find 
$x,y\in\p{2N+1}$ such that: \\
(a) $f^{T_{0}}=xy$ and $f^{T_{0}+\epsilon}=
x\cdot\phi_{i}^{N}(A_{i,j}^{\pm 1})\cdot y$ for some $j$;\\ 
(b) $x=\phi_{i}^{N}(x')$ for some $x'\in\p{2N}$.\\ Thus
\[ f^{T_{0}+\epsilon}=x\phi_{i}^{N}(A_{i,j}^{\pm 1})x^{-1}\cdot f^{T_{0}}=
 \phi_{i}^{N}(x'A_{i,j}^{\pm 1}{x'}^{-1})\cdot f^{T_{0}}.\]
Notice that conjugation by $x'$ maps $A_{i,j}$ to a product of $A_{i,j_{m}}$
for some set of $j_{m}$,  so $\phi_{i}^{N}(x'A_{i,j}^{\pm 1}{x'}^{-1})$ 
lies in $H^{T}$.


Similarly, if the rubber strip is  attached to the 
minimum, $f^{T_{0}}$ is multiplied on the right by some braid from $H^B$. 
In case two rubber strips intersect each other
we have to multiply by a product of two braids of such form; as above,
the product will lie in $H^T$ or $H^B$. (If one rubber strip is attached 
to a minimum and the other one to a maximum this product will automatically 
lie in the intersection $H^T\cap H^B$.)
Finally, when the isotopy is finished and all minima and maxima have arrived
back to their places what may happen is that some rubber strips may be twisted.
This corresponds to multiplications by some $A_{i,i+1}$ 
on the left  for $i$  even  and 
on the right for $i$  odd.

\bigskip

\noindent{\bf Acknowledgments.}
We would like to thank Mario Eudave for finding an important reference,
Sofia Lambropoulou and other organizers of the conference 
``Knots in Hellas '98'' who gave us a chance to meet, and Natig Atakishiev 
with whose pen a part of this paper was written. The second author was
partially supported by the Naval Academy Research Council.



{\small \begin{thebibliography}{99}

%
\bibitem[B1]{Bi1} J. Birman,
Braids, links, and mapping class groups, 
{\em Annals of Mathematics Studies,} No. 82. Princeton University Press, 
Princeton, N.J.; University of Tokyo Press, Tokyo, 1974.  
%
\bibitem[B2]{Bi2} J. Birman,
{\em On the stable equivalence of plat representations of knots and links}, 
Canad.\ J.\ Math.\ {\bf 28} (1976), no. 2, 264--290. 
%
\bibitem[Sch]{Sch} H. Schubert, 
{\em \"{U}ber eine numerische Knoteninvariante}, 
Math.\ Z.\ {\bf 61} (1954), 245--288.
%
\bibitem[St]{St} T. Stanford, 
{\em Vassiliev invariants and knots modulo pure braid subgroups},
math.GT/9805092.
%
\end{thebibliography}}
\end{document}